\documentclass[graybox]{svmult}

\usepackage{mathptmx}       
\usepackage{helvet}         
\usepackage{courier}        
\usepackage{type1cm}        
\usepackage{makeidx}         
\usepackage{graphicx}        
\usepackage{multicol}        
\usepackage[bottom]{footmisc}

\usepackage{amsfonts}
\usepackage{enumerate}
\usepackage{amsmath}
\usepackage{amssymb}
\usepackage{amsxtra}

\usepackage{graphicx, color}

\newcommand{\be}{\begin{equation}}
\newcommand{\ee}{\end{equation}}
\newcommand{\bea}{\begin{eqnarray}}
\newcommand{\eea}{\end{eqnarray}}
\newcommand{\beastar}{\begin{eqnarray*}}
\newcommand{\eeastar}{\end{eqnarray*}}

\def\E{{\mathbb E}}

\def\C{{\mathbb C}}
\def\R{{\mathbb R}}

\def\11{{\mathbb I}}

\def\CL{{\mathcal L}}

\def\delbar{{\overline \partial}}

\def\del{\partial}

\makeindex             

\begin{document}

\title*{Canonical connection on contact manifolds}
\author{Yong-Geun Oh and Rui Wang}
\institute{Yong-Geun Oh \at Center for Geometry and Physics, Institute for Basic Science (IBS),
77 Cheongam-ro, Nam-gu, Pohang-si, Gyeongsangbuk-do, Korea 790-784, \& POSTECH, Pohang, Korea
790-784,\\ \email{yongoh@ibs.re.kr}
\and Rui Wang \at Center for Geometry and Physics, Institute for Basic Science (IBS),
77 Cheongam-ro, Nam-gu, Pohang-si, Gyeongsangbuk-do, Korea 790-784, \email{rwang@ibs.re.kr}}
\maketitle

\abstract{We introduce a family of canonical affine connections on the contact manifold $(Q,\xi)$,
which is associated to each contact triad $(Q,\lambda,J)$ where
$\lambda$ is a contact form and $J:\xi \to \xi$ is an endomorphism with $J^2 = -id$
compatible to $d\lambda$. We call a particular one in this family the \emph{contact triad connection} of $(Q,\lambda,J)$
and prove its existence and uniqueness. The connection is canonical in that the pull-back connection
$\phi^*\nabla$ of a triad connection $\nabla$ becomes the triad connection of
the pull-back triad $(Q, \phi^*\lambda, \phi^*J)$ for any diffeomorphism $\phi:Q \to Q$. It also
preserves both the triad metric
$
g:= d\lambda(\cdot, J\cdot) + \lambda \otimes \lambda
$
and $J$ regarded as an endomorphism on $TQ = \mathbb{R}\{X_\lambda\}\oplus \xi$, and
is characterized by its torsion properties and the requirement that the contact
form $\lambda$ be holomorphic in the $CR$-sense. In particular, the connection
restricts to a Hermitian connection $\nabla^\pi$ on the Hermitian vector bundle
$(\xi,J,g_\xi)$ with $g_\xi = d\lambda(\cdot, J\cdot)|_{\xi}$, which
we call the \emph{contact Hermitian connection} of $(\xi,J,g_\xi)$.
These connections greatly simplify tensorial calculations
in the sequels \cite{oh-wang1}, \cite{oh-wang2} performed in the authors' analytic study of
the map $w$, called contact instantons, which satisfy the nonlinear elliptic system of equations
$$
\overline\partial^\pi w = 0, \, d(w^*\lambda \circ j) = 0
$$
of the contact triad $(Q,\lambda,J)$.
}

\section{Introduction}
\label{sec:intro}

Let $(Q,\xi)$ be a $2n+1$ dimensional contact manifold and a contact form $\lambda$
be given, which means that the contact distribution $\xi$ is given as $\ker \lambda$ and $\lambda\wedge(d\lambda)^n$
nowhere vanishes. On $Q$, the Reeb vector field $X_\lambda$ associated to the contact
form $\lambda$ is the unique vector field satisfying
$X_\lambda \rfloor \lambda = 1$ and $X_\lambda \rfloor d\lambda = 0$.
Therefore the tangent bundle $TQ$ has the splitting $TQ =\mathbb{R}\{X_\lambda\} \oplus \xi$.
We denote by
$
\pi: TQ \to \xi
$
the corresponding projection.

Now let $J$ be a complex structure on $\xi$, i.e., $J:\xi\to \xi$ with $J^2=-id|_\xi$.
We extend $J$ to $TQ$ by defining $J(X_\lambda)=0$.
We will use such $J: TQ \rightarrow TQ$ throughout the paper. Then we have $J^2 = -\Pi$ where
$\Pi: TQ \to TQ$ is the unique idempotent with $\operatorname{Im} \Pi = \xi$ and $\ker \Pi = \R\cdot X_\lambda$.
We note that we have the unique decomposition $h = \lambda(h) X_\lambda + \pi h$ for any $h \in TQ$ in terms of
the decomposition $TQ = \R \cdot X_\lambda \oplus \xi$.

\begin{definition}[Contact triad metric]
Let $(Q,\lambda,J)$ be a contact triad.
We call the metric defined by
$
g(h,k) : = \lambda(h)\lambda(k)
+ d\lambda(\pi h, J \pi k)
$ for any $h,\ k\in TQ$ the \emph{contact triad metric} associated to the triad $(Q,\lambda,J)$.
\end{definition}

The main purpose of the present paper is to introduce the notion of
the \emph{contact triad connection} of the triad $(Q,\lambda,J)$ which is
the contact analog to the Ehresman-Libermann's notion of \emph{canonical connection}
on the almost K\"ahler manifold $(M,\omega,J)$.
(See \cite{ehres-liber}, \cite{libermann54}, \cite{libermann55}, \cite{kobayashi}, \cite{gaudu} for general exposition
on the canonical connection.)
\begin{theorem}[Contact triad connection]\label{thm:mainI} Let $(Q,\lambda,J)$ be any contact triad of
contact manifold $(Q,\xi)$, and $g$ the contact triad connection.
Then there exists a unique affine connection $\nabla$ that has the
following properties:
\begin{enumerate}
\item $\nabla$ is a Riemannian connection of the triad metric.
\item The torsion tensor of $\nabla$ satisfies $T(X_\lambda, Y)=0$ for all $Y \in TQ$.
\item $\nabla_{X_\lambda} X_\lambda = 0$ and $\nabla_Y X_\lambda\in \xi$, for $Y\in \xi$.
\item $\nabla^\pi := \pi \nabla|_\xi$ defines a Hermitian connection of the vector bundle
$\xi \to Q$ with Hermitian structure $(d\lambda, J)$.
\item The $\xi$ projection, denoted by $T^\pi: = \pi T$, of the torsion $T$ has vanishing
$(1,1)$-component in its complexification, i.e., satisfies the following properties:
for all $Y$ tangent to $\xi$, $T^\pi(JY,Y) = 0$.
\item
For $Y\in \xi$,
$
\nabla_{JY} X_\lambda+J\nabla_Y X_\lambda = 0.
$
\end{enumerate}
We call $\nabla$ the contact triad connection.
\end{theorem}
Recall that the leaf space of Reeb foliations of the contact triad $(Q,\lambda,J)$
canonically carries a (non-Hausdorff) almost K\"ahler structure which we denote by
$
(\widehat Q,\widehat{d\lambda},\widehat J).
$
We would like to note that Axioms (4) and (5)
are nothing but properties of the canonical connection
on the tangent bundle of the (non-Hausdorff) almost
K\"ahler manifold $(\widehat Q, \widehat{d\lambda},\widehat J_\xi)$
lifted to $\xi$. (In fact, as in the almost K\"ahler case, vanishing of $(1,1)$-component
also implies vanishing of $(2,0)$-component and hence the torsion automatically
becomes $(0,2)$-type.)
On the other hand, Axioms (1), (2), (3)
indicate this connection behaves like the Levi-Civita connection when the Reeb direction $X_\lambda$ get involved.
Axiom (6) is an extra requirement to connect the information in $\xi$ part and $X_\lambda$ part,
which is used to dramatically simplify our calculation in \cite{oh-wang1}, \cite{oh-wang2}.

In fact, the contact triad connection is one of the $\R$-family of affine connections
satisfying Axioms (1) - (5) with (6) replaced by
$$
\nabla_{JY} X_\lambda + J \nabla_Y X_\lambda = c\ Y, \quad c \in \R.
$$
Contact triad connection corresponds to $c = 0$ and the connection $\nabla^{LC} +B_1$ (see Section \ref{sec:existence} for the expression of $B_1$)
corresponds to
$c = -1$. 

The contact triad connection (and also the whole $\R$-family) we construct here has naturality as stated below.
\begin{corollary}[Naturality] Let $\nabla$ be the contact triad connection of the triad $(Q,\lambda,J)$.
Then for any diffeomorphism $\phi: Q \to Q$, the pull-back connection $\phi^*\nabla$ is the contact triad
connection associated to the triad $(Q,\phi^*\lambda,\phi^*J)$.
\end{corollary}

While our introduction of Axiom (6) is motivated by our attempt to simplify the tensor
calculations \cite{oh-wang1}, it has a nice geometric interpretation in terms of
$CR$-geometry.  
(We refer to Definition \ref{defn:CR-holomorphic} for the definition
for $CR$-holomorphic $k$-forms.)

\begin{proposition} In the presence of other defining
properties of contact triad connection, Axiom (6) is equivalent to the statement that
$\lambda$ is holomorphic in the $CR$-sense.
\end{proposition}

\medskip
Some motivations of the study of the canonical connection are in order.
Hofer-Wysocki-Zehnder \cite{HWZ96,HWZ98} derived exponential decay estimates of
proper pseudoholomorphic curves with respect to the cylindrical almost complex structure associated to the
endomorphism $J: \xi \to \xi$ in symplectization
by bruit force coordinate calculations using some special
coordinates around the given Reeb orbit which is rather complicated.
Our attempt to improve the presentation of these decay estimates, using the tensorial language,
was the starting point of the research performed in the present paper.

We do this in  \cite{oh-wang1}, \cite{oh-wang2} by
considering a map $w:\dot \Sigma \to Q$ satisfying the equation
\be\label{eq:w}
\overline\partial^\pi w = 0, \, \quad d(w^*\lambda \circ j) = 0
\ee
\emph{without involving the function $a$ on the contact manifold $Q$
or the symplectization}. We call such a map a contact instanton.
We refer \cite{hofer-survey} for the origin of this equation in
contact geometry, as well as \cite{oh-wang1}, \cite{oh-wang2} for the detailed analytic study of priori
$W^{k,2}$-estimates and asymptotic convergence on punctured Riemann surfaces.

In the course of our studying the geometric analysis of such maps,
we need to simplify the tensorial calculations by choosing a special
connection as in the (almost) complex geometry. It turns out that for the purpose of
taking the derivatives of the map $w$ several times,
the contact triad connection on $Q$ is much more convenient and easier to keep track of
various terms arising from switching the order of derivatives
than the commonly used Levi-Civita connection. The advantage of the contact triad
connection will become even more apparent  in \cite{oh:sigmamodel} where the Fredholm theory
and the corresponding index computations in relation to the equations \eqref{eq:w} are developed.

There have been several literatures that studied special connections on contact manifolds, such as
\cite{tanno}, \cite{nicola}, \cite{stadtmuller}. We make some rough comparisons between these connections
and the contact triad connection introduced in this paper.

Although all the connections mentioned above are characterized by the torsion properties, one big difference
between ours and the ones in \cite{nicola}, \cite{stadtmuller} is that we don't require $\nabla J=0$,
but only $\nabla^\pi J=0$. Notice that $\nabla J=0$ is equivalent to both $\nabla^\pi J=0$
and $\nabla X_\lambda\in \R\cdot X_\lambda$. Together with the metric property, $\nabla J=0$
also implies $\nabla X_\lambda=0$, which is the requirement of the contact metric connection
studied in \cite[Def. 3.1]{nicola} as well as the so-called adapted connection considered
in \cite[Sec. 4]{stadtmuller}. Our contact triad connection doesn't satisfy this requirement
in general, and so is not in these families.

The connections considered in \cite{nicola}, \cite{stadtmuller} become the canonical connection
when lifted to the \emph{symplectization} as an almost K\"ahler manifold, while our connection and
the generalized Tanaka-Webster connection considered by Tanno \cite{tanno} are canonical
for the (non-Hausdorff) almost K\"ahler manifold $(\widehat Q, \widehat{d\lambda},\widehat J_\xi)$
lifted to $\xi$. (We remark that some other people named their connections the generalized
Tanaka-Webster connection with different meanings.)


Difference in our connection and Tanno's shows up in the torsion property of
$T(X_\lambda, \cdot)$ among others. It would be interesting to provide the classification of
the canonical connections in a bigger family that
includes both the contact triad connection and Tanno's generalized Tanaka-Webster connection.
Since the torsion of the triad connection is already reduced to the simplest one,
we expect that it satisfies better property on its curvature and get better results on the gauge
invariant studied in \cite{tanno}.

This paper is a simplified version of \cite{oh-wang},
to which we refer readers for the complete proofs of various results given in this paper.


\section{Review of the canonical connection of almost K\"ahler manifold}
\label{sec:can-Kahler}

We recall this construction of the canonical connection for
almost K\"ahler manifolds $(M,\omega,J)$.
A nice and exhaustive discussion on the general almost Hermitian connection is given by
Gauduchon in \cite{gaudu} to which we refer readers for more details.
(See also \cite{kobayashi}, \cite[Section 7.1]{oh:book}.)

Assume $(M,J,g)$ an almost Hermitian manifold, which means
$J$ is an almost complex structure $J$ and $g$ the metric satisfying
$
g(J\cdot ,J\cdot ) = g(\cdot, \cdot)
$.
An affine connection $\nabla$ is called $J$-linear if $\nabla J = 0$.
There always exists a $J$-linear connection for a given almost
complex manifold. We denote by $T$ the torsion tensor of $\nabla$.

\begin{definition}\label{defn:real-canonical} Let $(M,J,g)$ be an almost Hermitian manifold.
A $J$-linear connection is called a (the) canonical connection (or a (the) Chern connection) if
for any for any vector field $Y$ on $M$ there is
$T(JY,Y) = 0$.
\end{definition}

Recall that any $J$-linear connection extended to the complexification
$T_\C M = TM \otimes_{\mathbb{R}} \C$ complex linearly preserves the splitting
into $T^{(1,0)}M$ and $T^{(0,1)}M$.
Similarly we can extend the torsion tensor $T$ complex linearly which we denote by
$T_\C$.
Following the notation of \cite{kobayashi}, we denote
$
\Theta = \Pi' T_\C
$
the $T^{(1,0)}M$-valued two-form, where $\Pi'$ is the projection to $T^{(1,0)}M$.
We have the decomposition
$
\Theta = \Theta^{(2,0)} + \Theta^{(1,1)} + \Theta^{(0,2)}
$.
We can define the canonical connection in terms of the induced
connection on the complex vector bundle $T^{(1,0)}M \to M$.  The following lemma is easy to check by
definition.

\begin{lemma} An affine connection $\nabla$ on $M$ is a (the) canonical connection if and only if
the induced connection $\nabla$ on the complex vector bundle $T^{(1,0)}M$
has its complex torsion form $\Theta = \Pi' T_\C$ satisfy
$\Theta^{(1,1)} = 0$.
\end{lemma}

We particularly quote two theorems from Gauduchon \cite{gaudu}, Kobayashi \cite{kobayashi}.

\begin{theorem}\label{thm:ac-unique} On any almost Hermitian manifold $(M,J,g)$, there
exists a unique Hermitian connection $\nabla$ on $TM$ leading to the
canonical connection on $T^{(1,0)}M$. We call this connection the canonical Hermitian
connection of $(M,J,g)$.
\end{theorem}

We recall that $(M,J,g)$ is almost-K\"ahler if the fundamental two-form $\Phi = g(J \cdot, \cdot)$
is closed \cite{kobayashi-nomizu}.

\begin{theorem} Let $(M,J,g)$ be almost
K\"ahler and $\nabla$ be the canonical connection of $T^{(1,0)}M$. Then
$\Theta^{(2,0)} = 0$ in addition, and hence $\Theta$ is of type $(0,2)$.
\end{theorem}

\begin{remark}\label{rem:TJYY}It is easy to check by definition (or see \cite{gaudu}, \cite{kobayashi} for details) that $\Theta$ is of type $(0,2)$ is equivalent to say
that for all vector fields $Y, \, Z$ on $W$,
$T(JY,Z)= T(Y,JZ)$ and
$JT(JY,Z) =T(Y,Z)$.
\end{remark}


Now we describe one way of constructing the canonical connection
on an almost complex manifold described in \cite[Theorem 3.4]{kobayashi-nomizu}
which will be useful for our purpose of constructing the contact analog thereof
later. This connection has its torsion which satisfies
$
N = 4 T
$,
where $N$ is the Nijenhuis tensor of the almost complex structure $J$ defined as
$N(X, Y)= [JX, JY]-[X, Y]-J[X, JY]-J[JX, Y]$.
In particular, the complexification $\Theta = \Pi' T_\C$ is of $(0,2)$-type.

We now describe the construction of this canonical connection.
Let $\nabla^{LC}$ be the Levi-Civita connection.
Consider the standard averaged connection $\nabla^{av}$ of multiplication $J: TM \to TM$,
$$
\nabla^{av}_X Y : = \frac{\nabla^{LC}_X Y + J^{-1} \nabla^{LC}_X(JY)}{2}
=\nabla^{LC}_X Y - \frac{1}{2} J(\nabla^{LC}_XJ )Y.
$$
We then have the following Proposition stating that this connection becomes the canonical connection.
Its proof can be found in \cite[Theorem 3.4]{kobayashi-nomizu} or from section 2 \cite{gaudu}
with a little more strengthened argument by using \eqref{eq:sum-vanishing} for the metric property.

\begin{proposition}\label{prop:N=8T} Assume that $(M,g,J)$ is almost K\"ahler, i.e, the two-form
$\omega = g(J\cdot, \cdot)$ is closed. Then the average connection $\nabla^{av}$
defines the canonical connection of $(M,g,J)$, i.e., the connection
is $J$-linear, preserves the metric and its complexified torsion is
of $(0,2)$-type.
\end{proposition}

In fact, a more general construction of the canonical connection for almost Hermitian manifold
is given in \cite{kobayashi-nomizu}. We describe it and in later sections, we will give a contact analog of this construction.

Consider the tensor field $Q$ defined by
\be\label{eq:tensorQ}
4 Q(X,Y) = (\nabla^{LC}_{JY}J) X + J((\nabla^{LC}_Y J) X) + 2J((\nabla^{LC}_XJ)Y)
\ee
for vector fields $X, \, Y$ on $M$. It turns out that when $(M,g,J)$ is almost K\"ahler, i.e.,
the two form $g(J\cdot, \cdot)$ is closed, the sum of the first two terms vanish.
In general, $\nabla:=\nabla^{LC}-Q$ is the canonical connection of the almost Hermitian manifold.
In fact, we have the following lemma which explains the construction above for almost K\"ahler case.
\begin{lemma}[(2.2.10)\cite{gaudu}]\label{lem:sum-vanishing} Assume $(M,g,J)$ is almost K\"ahler. Then
\be\label{eq:sum-vanishing}
\nabla^{LC}_{JY}J + J(\nabla^{LC}_Y J) = 0
\ee
and so $Q(X,Y) = \frac{1}{2} J((\nabla^{LC}_XJ )Y)$.
\end{lemma}

\section{Definition of the contact triad connection and its consequences}
\label{sec:can-connection}



In this section, we associate a particular type of affine connection on $Q$
to the given contact triad
$
(Q,\lambda,J)
$
which we call \emph{the contact triad connection} of the triple.


We recall
$
TQ = \R \{X_\lambda\} \oplus \xi$,
and denote by $\pi: TQ \to \xi$ the projection. Under this splitting, we
may regard a section $Y$ of $\xi \to Q$ as a vector field $Y \oplus 0$. We will
just denote the latter by $Y$ with slight abuse of notation.
Define $\nabla^\pi$ the connection of the bundle $\xi\rightarrow Q$ by
$\nabla^\pi Y=\pi \nabla Y$.

\begin{definition}[Contact triad connection]\label{defn:can-connection} We call an affine connection $\nabla$ on $Q$
the \emph{contact triad connection} of the contact triad $(Q,\lambda,J)$, if
it satisfies the following properties:
\begin{enumerate}
\item $\nabla^\pi$ is a Hermitian connection of the Hermitian bundle $\xi$ over the contact manifold $Q$ with Hermitian structure
$(d\lambda, J)$.
\item The $\xi$ projection, denoted by $T^\pi: = \pi T$, of the torsion $T$ satisfies the following properties:
for all $Y$ tangent to $\xi$,
$
T^\pi(JY,Y) = 0$.
\item $T(X_\lambda, Y)=0$ for all $Y \in TQ$.
\item $\nabla_{X_\lambda} X_\lambda = 0$ and  $\nabla_Y X_\lambda\in \xi$, for $Y\in \xi$.
\item For $Y\in \xi$,
$
\nabla_{JY} X_\lambda+J\nabla_Y X_\lambda= 0.
$
\item For any $Y, \ Z\in \xi$,
$
\langle \nabla_Y X_\lambda, Z\rangle + \langle X_\lambda, \nabla_Y Z\rangle =0.
$
\end{enumerate}
\end{definition}
It follows from the definition that
the contact triad connection is a Riemannian connection of the triad metric.
(The statements of this definition are equivalent to those given in the introduction. We
state properties of contact triad connection here as above which are
organized in the way how they are used in the proofs of uniqueness and existence.)
%

By the second part of Axiom (4), the covariant derivative $\nabla X_\lambda$ restricted to $\xi$
can be decomposed into
$
\nabla X_\lambda = \del^{\nabla}X_\lambda + \delbar^{\nabla}X_\lambda
$,
where $\del^{\nabla}X_\lambda$ (respectively, $\delbar^{\nabla}X_\lambda$) is $J$-linear (respectively,
$J$-anti-linear part). Axiom (6) then is nothing but requiring that $\del^{\nabla}X_\lambda = 0$, i.e.,
$X_\lambda$ is anti $J$-holomorphic in the $CR$-sense. (It appears that this explains the reason why Axiom (5)
gives rise to dramatic simplification in our tensor calculations performed in \cite{oh-wang1}.)

One can also consider similar decompositions
of one-form $\lambda$. For this, we need some digression.
Define $J\alpha$ for a $k$-form $\alpha$ by the formula
$
J\alpha(Y_1, \cdots, Y_k) = \alpha(JY_1, \cdots , JY_k)
$.
\begin{definition}\label{defn:CR-holomorphic} Let $(Q,\lambda,J)$ be a contact triad. We call a $k$-form is $CR$-holomorphic
if $\alpha$ satisfies
\bea
\nabla_{X_\lambda} \alpha & = & 0, \label{eq:Xlambdaalpha}\\
\nabla_Y \alpha + J \nabla_{JY} \alpha & = & 0 \quad \text{for }\, Y \in \xi.\label{eq:Yalpha}
\eea
\end{definition}

\begin{proposition}\label{prop:CR} Axiom (5)  is equivalent to the statement that
$\lambda$ is holomorphic in the $CR$-sense in the presence of other defining
properties of contact triad connection.
\end{proposition}
\begin{proof} We first prove $\nabla_{X_\lambda} \lambda  = 0$ by evaluating it
against vector fields on $Q$. For $X_\lambda$, the first half of Axiom (4) gives rise to
$
\nabla_{X_\lambda}\lambda (X_\lambda) = - \lambda(\nabla_{X_\lambda}X_\lambda) = 0
$.
For the vector field $Y \in \xi$, we compute
\beastar
\nabla_{X_\lambda} \lambda(Y) & = & - \lambda(\nabla_{X_\lambda}Y)\\
& =& -\lambda(\nabla_Y X_\lambda + [X_\lambda,Y] + T(X_\lambda,Y)) \\
& = &  -\lambda(\nabla_Y X_\lambda) - \lambda([X_\lambda,Y]) - \lambda(T(X_\lambda,Y)).
\eeastar
Here the third term vanishes by Axiom (3), the first term by the second part of Axiom (4) and the second term
vanishes since
$$
\lambda([X_\lambda,Y]) = \lambda(\CL_{X_\lambda}Y) = X_\lambda[\lambda(Y)] - \CL_{X_\lambda}\lambda (Y)
= 0 - 0 = 0.
$$
Here the first vanishes since $Y \in \xi$ and the second because $ \CL_{X_\lambda}\lambda=0$ by
the definition of the Reeb vector field. This proves \eqref{eq:Xlambdaalpha}.

We next compute $J\nabla_Y \lambda$ for $Y \in \xi$. For a vector field $Z \in \xi$,
\beastar
(J\nabla_Y\lambda)(Z)&=&(\nabla_Y \lambda)(JZ)
= \nabla_Y(\lambda(JZ))-\lambda(\nabla_Y(JZ))
= -\lambda(\nabla_Y(JZ))
\eeastar
since $\lambda(JZ) = 0$ for the last equality. Then by the definitions of the Reeb vector field and the triad
metric and the skew-symmetry of $J$, we derive
\beastar
-\lambda(\nabla_Y(JZ))
= -\langle \nabla_Y(JZ), X_\lambda\rangle
= \langle JZ, \nabla_Y X_\lambda\rangle
= -\langle Z, J\nabla_Y X_\lambda\rangle.
\eeastar
Finally, applying (6), we obtain
\beastar
-\langle Z, J\nabla_Y X_\lambda\rangle
= \langle Z, \nabla_{JY} X_\lambda\rangle
= -\langle \nabla_{JY}Z,  X_\lambda\rangle
= -\lambda(\nabla_{JY}Z)
= (\nabla_{JY}\lambda)(Z).
\eeastar
Combining the above, we have derived
$
J(\nabla_Y\lambda)(Z)=\nabla_{JY}\lambda(Z)
$
for all $Z \in \xi$. On the other hand, for $X_\lambda$, we evaluate
$$
J(\nabla_Y\lambda)(X_\lambda) = \nabla_Y\lambda(JX_\lambda)  = \nabla_Y\lambda(0) = 0.
$$
We also compute
$
\nabla_{JY}\lambda(X_\lambda) = \CL_{JY}(\lambda(X_\lambda)) -\lambda(\nabla_{JY} X_\lambda)).
$
The first term vanishes since $\lambda(X_\lambda) \equiv 1$ and the second
vanishes since $\nabla_{JY} X_\lambda \in \xi$ by the second part of Axiom (4).
Therefore we have derived \eqref{eq:Yalpha}.

 Combining \eqref{eq:Xlambdaalpha} and
\eqref{eq:Yalpha}, we have proved that
Axiom (5) implies $\lambda$ is holomorphic in the $CR$-sense. The converse can be
proved by reading the above proof backwards.
\end{proof}
From now on, when we refer Axioms, we mean the properties in Definition \ref{defn:can-connection}.
One very interesting consequence of this uniqueness is the following naturality result of
the contact-triad connection.

\begin{theorem}[Naturality]\label{thm:naturality}
Let $\nabla$ be the contact triad connection of the triad $(Q,\lambda,J)$.
For any diffeomorphism $\phi: Q \to Q$,
the pull-back connection $\phi^*\nabla$ is the contact triad
connection associated to the triad $(Q,\phi^*\lambda,\phi^*J)$.
\end{theorem}
\begin{proof} A straightforward computation shows that the pull-back connection $\phi^*\nabla$
satisfies all Axioms $(1) - (6)$  for the triad $(Q, \phi^*\lambda, \phi^*J)$.
Therefore by the uniqueness, $\phi^*\nabla$ is the canonical connection.

\end{proof}
\begin{remark}\label{rem:naturaity}
An easy examination of the proof of Theorem \ref{thm:naturality}
shows that the naturality property stated in Theorem \ref{thm:naturality}
also holds for the one-parameter family of connections for all $c \in \R$ (see Section \ref{sec:uniqueness})
among which the canonical connection corresponds to $c = 0$.
\end{remark}

\section{Proof of the uniqueness of the contact triad connection}
\label{sec:uniqueness}

In this section, we give the uniqueness proof by analyzing the first structure equation and
showing how every axiom determines the connection one forms. In the next two sections, we explicitly
construct a connection by carefully examining properties of the
Levi-Civita connection and modifying the constructions in \cite{kobayashi-nomizu},
\cite{kobayashi} for  the canonical connection, and then show it satisfies all the requirements
and thus the unique contact triad connection.

We are going to prove the existence and uniqueness for a more general family of connections.
First, we generalize the Axiom (5) to the following Axiom:
For $Y\in \xi$,
\be\label{eq:JYJnablainR}
\nabla_{JY}X_\lambda+J\nabla_Y X_\lambda\in \mathbb{R}\cdot Y,
\ee
and we denote by Axiom (5; $c$): For a given $c \in \mathbb{R}$,
\be\label{eq:JYJnabla-cY}
\nabla_{JY}X_\lambda+J\nabla_Y X_\lambda = c\, Y, \quad Y \in \xi.
\ee
In particular, Axiom (5) corresponds to Axiom (5; $0$).

\begin{theorem}\label{thm:cont-unique}
For any $c \in \mathbb{R}$, there exists a unique connection satisfies
Axiom (1)-(4), (6) and (5; $c$).
\end{theorem}

\begin{proof}(Uniqueness)

Choose a moving frame of $TQ = \R\{X_\lambda\} \oplus \xi$ given by
$
\{X_\lambda, E_1,\cdots, E_n, J E_1,\cdots, J E_n\}
$
and denote its dual co-frame by
$
\{\lambda, \alpha^1,\cdots, \alpha^n, \beta^1,\cdots, \beta^n\}
$.
(We use the Einstein summation convention to denote the sum of upper indices and lower indices in this paper.)
Assume the connection matrix is $(\Omega^i_j)$, $i, j=0, 1,..., 2n$, and we write the first structure equations as follows
\begin{eqnarray*}
d\lambda&=&-\Omega^0_0\wedge \lambda-\Omega^0_k\wedge\alpha^k-\Omega^0_{n+k}\wedge\beta^k+T^0\\
d\alpha^j&=&-\Omega^j_0\wedge \lambda-\Omega^j_k\wedge\alpha^k-\Omega^j_{n+k}\wedge\beta^k+T^j\\
d\beta^j&=&-\Omega^{n+j}_0\wedge \lambda-\Omega^{n+j}_k\wedge\alpha^k-\Omega^{n+j}_{n+k}\wedge\beta^k+T^{n+j}
\end{eqnarray*}
Throughout the section, if not stated otherwise, we let  $i$, $j$ and $k$ take values from $1$ to $n$.
Denote
$$
\Omega^u_v=\Gamma^u_{0, v}\lambda+\Gamma^u_{k, v}\alpha^k+\Gamma^u_{n+k, v}\beta^k
$$
where $u, v=0, 1, \cdots, 2n$.
We will analyze each axiom in Definition \ref{defn:can-connection} and show how they set down the matrix of connection one forms.

We first state that Axioms (1) and (2) uniquely determine $(\Omega^i_j |_\xi)_{i,j=1,\cdots,2n}$.
This is exactly the same as Kobayashi's proof for the uniqueness of Hermitian connection
given in \cite{kobayashi}.
To be more specific, we can restrict the first structure equation to $\xi$ and get the following equations
for $\alpha$ and $\beta$ since $\xi$ is the kernel of $\lambda$.

\beastar
d\alpha^j&=&-\Omega^j_k|_\xi\wedge\alpha^k-\Omega^j_{n+k}|_\xi\wedge\beta^k+T^j|_\xi\\
d\beta^j&=&-\Omega^{n+j}_k|_\xi\wedge\alpha^k-\Omega^{n+j}_{n+k}|_\xi\wedge\beta^k+T^{n+j}|_\xi
\eeastar

We can see $(\Omega^i_j |_\xi)_{i,j=1,\cdots,2n}$ is skew-Hermitian from Axiom (1).  We also notice that
from the Remark \ref{rem:TJYY} that Axiom (2)
is equivalent to say that $\Theta^{(1,1)}=0$, where $\Theta= \Pi' T_\mathbb{C}$. Then one
can strictly follow Kobayashi's proof of Theorem \ref{thm:ac-unique} in \cite{kobayashi} and get $(\Omega^i_j |_\xi)_{i,j=1,\cdots,2n}$
are uniquely determined. For this part, we refer readers to the proofs of \cite[Theorem 1.1 and 2.1]{kobayashi}.

In the rest of the proof, we will clarify how the Axioms (3), (4), (5;$c$), (6) uniquely determine $\Omega^0_{\cdot}$, $\Omega^{\cdot}_0$ and $(\Omega^i_j (X_\lambda))_{i,j=1,\cdots,2n}$.
Compute the first equality in Axiom (4) and we get
$$\nabla_{X_\lambda}{X_\lambda}=\Gamma^0_{0, 0}X_{\lambda}+\Gamma^k_{0, 0}E_k+\Gamma^{n+k}_{0, 0}JE_k=0.$$
Hence
\be
\Gamma^0_{0, 0}=0, \quad \Gamma^k_{0, 0}=0, \quad \Gamma^{n+k}_{0, 0}=0\label{eq:gamma-000}
\ee
The second claim in Axiom (4) is equal to say
\bea
\nabla_{E_k}X_{\lambda}\in \xi,\quad
\nabla_{JE_k}X_{\lambda}\in \xi.
\eea
Similar calculation shows that
\be
\Gamma^0_{k, 0}=0, \quad
\Gamma^0_{n+k, 0}=0\label{eq:gamma-0k0}.
\ee
Now the first vanishing in \eqref{eq:gamma-000} together with \eqref{eq:gamma-0k0} uniquely settle down
$$\Omega^0_0=\Gamma^0_{0, 0}\lambda+\Gamma^0_{k, 0}\alpha^k+\Gamma^0_{n+k, 0}\beta^k=0.$$
The vanishing of second and third equality in \eqref{eq:gamma-000} will be used to determine $\Omega_0$ in the later part.
From Axiom (3), we can get
\bea
\Gamma^k_{j, 0}-\Gamma^k_{0,j}&=&\langle [E_j, X_\lambda], E_k\rangle=-\langle \CL_{X_\lambda}E_j, E_k\rangle\label{eq:condition-31}\\
\Gamma^k_{n+j, 0}-\Gamma^k_{0,n+j}&=&\langle [JE_j, X_\lambda], E_k\rangle=-\langle \CL_{X_\lambda}(JE_j), E_k\rangle\label{eq:condition-32}
\eea
and
\bea
\Gamma^{n+k}_{j, 0}-\Gamma^{n+k}_{0,j}&=&\langle [E_j, X_\lambda], JE_k\rangle=-\langle \CL_{X_\lambda}E_j, JE_k\rangle\label{eq:condition-33}\\
\Gamma^{n+k}_{n+j, 0}-\Gamma^{n+k}_{0,n+j}&=&\langle [E_j, X_\lambda], JE_k\rangle=-\langle \CL_{X_\lambda}(JE_j), JE_k\rangle\label{eq:condition-34}.
\eea
From Axiom (5; $c$), we have
\bea
\Gamma^k_{j, 0}+\Gamma^{n+k}_{n+j, 0}&=&0\label{eq:condition-51}\\
\Gamma^{n+k}_{j, 0}-\Gamma^k_{n+j, 0}&=&-c\delta_{j,k}\label{eq:condition-52}.
\eea

Now we show how to determine $\Omega^j_0$ for $j=1, \dots, 2n$. For this purpose, we calculate $\Gamma^k_{j, 0}$.
First, by using \eqref{eq:condition-51}, we write
$\Gamma^k_{j, 0}=\frac{1}{2}\Gamma^k_{j, 0}-\frac{1}{2}\Gamma^{n+k}_{n+j, 0}$.

Furthermore, using \eqref{eq:condition-31} and \eqref{eq:condition-34} , we have
\beastar
\Gamma^k_{j, 0}&=&\frac{1}{2}\Gamma^k_{j, 0}-\frac{1}{2}\Gamma^{n+k}_{n+j, 0}\\
&=&\frac{1}{2}(\Gamma^k_{0,j}-\langle \CL_{X_\lambda}E_j, E_k\rangle)-\frac{1}{2}(\Gamma^{n+k}_{0, n+j}-\langle \CL_{X_\lambda}(JE_j), JE_k\rangle)\\
&=&\frac{1}{2}(\Gamma^k_{0,j}-\Gamma^{n+k}_{0,n+j})-\frac{1}{2}(\langle \CL_{X_\lambda}E_j, E_k\rangle-\langle \CL_{X_\lambda}(JE_j), JE_k\rangle)\\
&=&\frac{1}{2}(\Gamma^k_{0,j}-\Gamma^{n+k}_{0,n+j})-\frac{1}{2}\langle \CL_{X_\lambda}E_j+J\CL_{X_\lambda}(JE_j), E_k\rangle\\
&=&\frac{1}{2}(\Gamma^k_{0,j}-\Gamma^{n+k}_{0,n+j})-\frac{1}{2}\langle J(\CL_{X_\lambda}J)E_j, E_k\rangle\\
&=&\frac{1}{2}(\Gamma^k_{0,j}-\Gamma^{n+k}_{0,n+j})+\frac{1}{2}\langle (\CL_{X_\lambda}J)JE_j, E_k\rangle
\eeastar
Notice the first term vanishes by Axiom (2). In particular, that is from $\nabla_{X_\lambda}J=0$.
Hence we get
\bea
\Gamma^k_{j, 0}=\frac{1}{2}\langle (\CL_{X_\lambda}J)JE_j, E_k\rangle.
\eea

Following the same idea, we use \eqref{eq:condition-52} and will get
\beastar
\Gamma^{n+k}_{j, 0}
=-\frac{1}{2}c\delta_{jk}+\frac{1}{2}\langle (\CL_{X_\lambda}J)JE_j, JE_k\rangle.
\eeastar
Then substituting this into \eqref{eq:condition-51} and \eqref{eq:condition-52}, we have
$$
\Gamma^{k}_{n+j, 0} = \frac{1}{2}c\delta_{jk}+\frac{1}{2}\langle (\CL_{X_\lambda}J)JE_j, JE_k\rangle
= \frac{1}{2}c\delta_{jk}-\frac{1}{2}\langle (\CL_{X_\lambda}J)E_j, E_k\rangle.
$$
and
$$
\Gamma^{n+k}_{n+j, 0}= -\frac{1}{2}\langle (\CL_{X_\lambda}J)JE_j, E_k\rangle
= \frac{1}{2}\langle (\CL_{X_\lambda}J)E_j, JE_k\rangle.
$$

Together with \eqref{eq:gamma-000},
$\Omega_0$ is uniquely determined by this way.

Furthermore \eqref{eq:condition-31},\eqref{eq:condition-32},\eqref{eq:condition-33}
and \eqref{eq:condition-34}, uniquely determine $\Omega^i_j(X_\lambda)$ for $i,j=1, \dots, 2n$.

Notice that for any $Y\in \xi$, we derive
$
\nabla_{X_\lambda}Y \in \xi
$
from Axiom (3). This is because the axiom implies $\nabla_{X_\lambda}Y =\nabla_Y X_\lambda+\CL_{X_\lambda}Y$
and the latter is contained in $\xi$: the second part of Axiom (4) implies $\nabla_Y X_\lambda \in \xi$ and
the Lie derivative along the Reeb vector field  preserves
the contact structure $\xi$. It then follows that $\Gamma^0_{0, l}=0$ for $l=1, \dots, 2n$.
At the same time, Axiom (6) implies
$
\Gamma^0_{j,k}=-\Gamma^k_{j,0}.
$
for $j,\ k=1, \dots, 2n$. Hence together with \eqref{eq:gamma-0k0}, $\Omega^0$ is uniquely determined.
This finishes the proof.
\end{proof}

We end this section by giving a summary of the procedure we take in the proof
of uniqueness which actually indicates a way how to construct this connection in later sections.

First, we use the Hermitian connection property, i.e., Axiom (1) and torsion property
Axiom (2), i.e., $T^\pi|_\xi$ has vanishing $(1,1)$ part, to uniquely fix the connection on $\xi$ projection of $\nabla$ when taking values on $\xi$.

Then we use the metric property
$
\langle X_\lambda, \nabla_Y Z\rangle +\langle \nabla_Y X_\lambda, Z\rangle=0,
$
for any $Y, Z\in \xi$, to determine the $X_\lambda$ component of $\nabla$ when taking values in $\xi$.

To do this, we need the information of $\nabla_Y X_\lambda$. As mentioned before the second part of
Axiom (4) enables us to decompose
$
\nabla X_\lambda = \del^{\nabla} X_\lambda + \delbar^{\nabla}X_\lambda
$.
The requirement $\nabla_{X_\lambda}J=0$ in Axiom (1) implies
$
\nabla_{X_\lambda}(JY)-J\nabla_{X_\lambda}Y=0.
$
Axiom (3), the torsion property $T(X_\lambda, Y)=0$, then interprets this one into
$$
\nabla_{JY}X_\lambda-J\nabla_Y X_\lambda=-(\CL_{X_\lambda}J)Y
$$
which is also equivalent to saying
\be\label{eq:delbarnablaX}
J \delbar^{\nabla}_Y X_\lambda = \frac{1}{2}(\CL_{X_\lambda}J)Y \quad \text{or } \, \delbar^{\nabla}_Y X_\lambda = \frac{1}{2}(\CL_{X_\lambda}J)J Y.
\ee
It turns out that we can vary Axiom (5) by replacing it to
(5;$c$)
\be\label{eq:delnablaX}
\nabla_{JY}X_\lambda+J\nabla_Y X_\lambda=cY, \quad\text{or equivalently } \, \del^{\nabla}_Y X_Y = \frac{c}2 Y
\ee
for any given real number $c$. This way we shall have one-parameter family of
affine connections parameterized by $\R$ each of which satisfies Axioms (1) - (4) and (6)
with (5) replaced by (5;$c$).

When $c$ is fixed, i.e., under Axiom (5; $c$),
we can uniquely determine $\nabla_Y X_\lambda$ to be
$$
\nabla_Y X_\lambda=-\frac{1}{2}cJY+\frac{1}{2}(\CL_{X_\lambda}J)JY.
$$
Therefore, $\nabla_Y, \, Y \in \xi$ is uniquely determined in this process by getting the formula of $\nabla_Y X_\lambda$
when combined with the torsion property.
Then the remaining property $\nabla_{X_\lambda}X_\lambda=0$ now completely determines the connection.

\section{Properties of the Levi-Civita Connection on contact manifolds}
\label{sec:LC}

From the discussion in previous sections, the only thing left to do for the existence
of the contact triad connection is to globally define a connection such that
it can patch the $\xi$ part of $\nabla|_\xi$ and the $X_\lambda$ part of it.
In particular, we seek for a connection that satisfies the following properties:
\begin{enumerate}
\item it satisfies all the algebraic properties of
the canonical connection of almost K\"ahler manifold \cite{kobayashi} when restricted to $\xi$.
\item it satisfies metric property and has vanishing torsion in $X_\lambda$ direction.
\end{enumerate}
%
The presence of such a construction is a manifestation of delicate interplay
between the geometric structures $\xi$, $\lambda$, and $J$
in the geometry of contact triads $(Q,\lambda,J)$. In this regard, the closeness of $d\lambda$
and the definition of Reeb vector field $X_\lambda$
play important roles. In particular $d\lambda$ plays the role similar to that of the fundamental
two-form $\Phi$ in the case of almost K\"ahler manifold \cite{kobayashi-nomizu} (in a non-strict sense)
in that it is closed.

This interplay is reflected already in several basic
properties of the Levi-Civita connection of the contact triad metric exposed in this section.
We list these properties but skip most proofs of them in this section since most results are well-known in Blair's book
\cite{blair}. We also refer readers to \cite{oh-wang} for the complete proof with the same convention.


Recall that we have extend $J$
to $TQ$ by defining $J(X_\lambda)=0$. Denote by $\Pi: TQ \to TQ$ the
idempotent associated to the projection $\pi:TQ \to \xi$, i.e., the endomorphism satisfying
$\Pi^2 = \Pi$, $\operatorname{Im} \Pi = \xi$, and $\ker \Pi = \R\{X_\lambda\}$.

We have now
$J^2=-\Pi$. Moreover,  for any connection $\nabla$ on $Q$,
\bea
 (\nabla J)J=-(\nabla \Pi) -J(\nabla J)\label{eq:J1}.
\eea
Notice for $Y\in \xi$, we have
\be
\Pi(\nabla \Pi)Y=0, \quad
(\nabla \Pi)X_\lambda = -\Pi \nabla X_\lambda.\label{eq:J2}
\ee

Denote the triad metric $g$ as $\langle \cdot, \cdot\rangle$.
By definition, we have
\beastar
\langle X, Y\rangle &=& d\lambda (X, JY)+\lambda(X)\lambda(Y)\\
d\lambda(X, Y)&=&d\lambda(JX, JY)
\eeastar
which gives rise to the following identities

\begin{lemma}\label{lemma:contact-J-compatible} For all $X, \, Y$ in $TQ$,
$
\langle JX, JY\rangle = d\lambda(X, JY)$,
$\langle X, JY\rangle = -d\lambda(X, Y)$, and
$\langle JX, Y\rangle =  -\langle X, JY\rangle
$.
\end{lemma}

However, we remark
$
\langle JX, JY\rangle \neq \langle X, Y\rangle
$
in general now, and hence there is no obvious analog of the fundamental 2 form $\Phi$
defined as in \cite{kobayashi-nomizu}
for the contact case. This is the main reason that is responsible
for the differences arising in the various relevant formulae
 between the contact case and the almost Hermitian case.

The following preparation lemma says that the linear operator $\CL_{X_\lambda}J$ is
symmetric with respect to the metric $g=\langle \cdot, \cdot \rangle$. 
\begin{lemma}[Lemma 6.2 \cite{blair}]\label{lem:Lie-symmetry}
For $Y, Z \in \xi$,
$\langle (\CL_{X_\lambda}J)Y, Z\rangle=\langle Y, (\CL_{X_\lambda}J)Z\rangle$.
\end{lemma}

The following simple but interesting lemma shows that the Reeb foliation is
a geodesic foliation for the Levi-Civita connection (and so for
the contact triad connection) of the contact triad metric.

\begin{lemma}[\cite{blair}]\label{lem:LC-XX} For any vector field $Z$ on $Q$,
\be\label{eq:nablaZX}
\nabla^{LC}_Z X_\lambda \in \xi,
\ee
and
\be\label{eq:Xgeodesic}
\nabla^{LC}_{X_\lambda}X_\lambda = 0.
\ee
\end{lemma}

Next we state the following lemma which is the contact analog to the Prop 4.2 in \cite{kobayashi-nomizu}
for the almost Hermitian case. The proof of this lemma can be also extracted from
\cite[Corollary 6.1]{blair} and so we skip it but refer \cite{oh-wang} for details.

\begin{lemma}\label{lem:inner-product-N-contact} Consider the Nijenhuis tensor $N$ defined by
\beastar
N(X, Y)= [JX, JY]-[X, Y]-J[X, JY]-J[JX, Y]
\eeastar
as in the almost complex case. For all $X, \, Y$ and $Z$ in $TQ$,
\beastar
2\langle (\nabla^{LC}_X J)Y, Z \rangle&=& \langle N(Y, Z), JX\rangle \\
&{}&-\langle JX, JY\rangle \lambda(Z)+\langle JX, JZ\rangle \lambda(Y)
\eeastar
\end{lemma}

In particular, we obtain the following corollary.
\begin{corollary}\label{cor:nablaLCJ}
For $Y, Z\in \xi$,
\beastar
2\langle (\nabla^{LC}_{Y} J)X_\lambda, Z \rangle&=&-\langle (\CL_{X_\lambda}J)Z, Y\rangle +\langle Y, Z\rangle\\
2\langle (\nabla^{LC}_{Y} J)Z, X_\lambda \rangle &=& \langle (\CL_{X_\lambda}J)Z, Y\rangle -\langle Y, Z\rangle\\
2\langle (\nabla^{LC}_{X} J)Y, Z \rangle&=& \langle N(Y, Z), JX\rangle .
\eeastar
\end{corollary}
\begin{proof}
This is a direct corollary from Lemma \ref{lem:inner-product-N-contact} except that we also use
\bea
N(X_\lambda, Z)&=&-J(\CL_{X_\lambda}J)Z\label{eq:NX1}\\
N(Z, X_\lambda)&=&J(\CL_{X_\lambda}J)Z.\label{eq:NX2}
\eea
for the first two conclusions.
\end{proof}

Straightforward calculations give the following lemma which is the contact analog of the fact that the
Nijenhuis tensor is of $(0,2)$-type. 
\begin{lemma}\label{lem:N-property-contact}
For $Y, Z\in \xi$,
\beastar
JN(Y, JZ)-\Pi N(Y, Z)&=&0\\
\Pi N(Y, JZ)+\Pi N(Z, JY)&=&0.
\eeastar
\end{lemma}

Together with the last equality in Corollary \ref{cor:nablaLCJ} and Lemma \ref{lem:N-property-contact},
we obtain the following lemma, which  is the contact analog to Lemma \ref{lem:sum-vanishing}.

\begin{lemma}
\be\label{eq:Q-vanishing}
\Pi(\nabla^{LC}_{JY}J)X+J(\nabla^{LC}_Y J)X=0.
\ee
\end{lemma}

The following result is an immediate but important corollary of Corollary \ref{cor:nablaLCJ} and the property
$\nabla_{X_\lambda}X_\lambda = 0$ of $X_\lambda$, which plays an essential role in our construction of the
contact triad connection. 

\begin{proposition}[Corollary 6.1 \cite{blair}]\label{prop:LCXJ}
$
\nabla^{LC}_{X_\lambda} J=0
$.
\end{proposition}

The following is equivalent to the second part of Lemma 6.2 \cite{blair} after taking
into consideration of different sign convention of the definition of compatibility of
$J$ and $d\lambda$. 

\begin{lemma}[Lemma 6.2 \cite{blair}]\label{lem:YXlambda-LC}
For any $Y\in \xi$, we have
$
\nabla^{LC}_Y X_\lambda=\frac{1}{2}JY+\frac{1}{2}(\CL_{X_\lambda}J)JY$.
\end{lemma}

\section{Existence of the contact triad connection}
\label{sec:existence}

In this section, we establish the existence theorem of the contact triad connection in two stages.


Before we give the construction, we first remark the relationship between the connections of two different $c$'s.
Denote by $\nabla^{\lambda;c}$ the unique connection associated to the constant $c$, which we are going to construct. The following proposition shows that $\nabla^{\lambda;c}$ and $\nabla^{\lambda;c'}$ for two different nonzero constants
with the same parity are essentially the same in that it arises from the scale change of the contact form.
We skip the proof since it is straightforward.

\begin{proposition} Let $(Q,\lambda,J)$ be a contact triad and consider the triad
$(Q,a \lambda, J)$ for a constant $a > 0$. Then
$
\nabla^{a\lambda;1} = \nabla^{\lambda;a}
$.
\end{proposition}

In regard to this proposition, one could say that for each given
contact structure $(Q,\xi)$, there are essentially two inequivalent $\nabla^0, \, \nabla^1$
(respectively three, $\nabla^0, \, \nabla^1$ and $\nabla^{-1}$, if one fixes the orientation) choice of triad connections
for each given projective equivalence class of the contact triad $(Q,\lambda,J)$.
In this regard, the connection $\nabla^0$ is
essentially different from others in that this argument of scaling procedure of contact form $\lambda$
does not apply to the case $a = 0$ since it would lead to the zero form $0\cdot \lambda$.
This proposition also reduces the construction
essentially two connections of $\nabla^{\lambda;0}$ and
$\nabla^{\lambda;1}$ (or $\nabla^{\lambda;-1}$).

In the rest of this section, we will explicitly construct $\nabla^{\lambda;-1}$ and $\nabla^{\lambda;c}$
in two stages, by construct the potential tensor $B$ from the Levi-Civita connection,
i.e., by adding suitable tensors $B$ to get $\nabla^B=\nabla^{LC}+B:=\nabla^{LC}+B_1+B_2$.

In the first stage, 
motivated by the construction of the canonical connection on almost K\"ahler manifold and use
the properties of the Levi-Civita connection we extracted in the previous section,
we construct the connection $\nabla^{tmp;1}$
and show that it satisfies Axioms (1)-(4), (5;$-1$), (6).

In the second  stage, we modify $\nabla^{tmp;1}$ to get $\nabla^{tmp;2}$
by deforming the property (5;$-1$) thereof to (5;$c$) leaving other properties of $\nabla^{tmp;1}$ intact.
This $\nabla^{tmp;2}$ then satisfies all the axioms in Definition \ref{defn:can-connection}.
%

\subsection{Modification 1; $\nabla^{tmp;1}$}
\label{subsec:modify1}
\par


Define an affine connection $\nabla^{tmp;1}$ by the formula
$$
\nabla^{tmp;1}_{Z_1}Z_2 = \nabla^{LC}_{Z_1} Z_2 -  \Pi P(\Pi Z_1,  \Pi Z_2)
$$
where the bilinear map $P: \Gamma(TQ) \times \Gamma(TQ) \to \Gamma(TQ)$
over $C^\infty(Q)$ is defined by
\be\label{eq:tensorP}
4 P(X,Y) = (\nabla^{LC}_{JY}J) X + J((\nabla^{LC}_Y J) X) + 2J((\nabla^{LC}_XJ)Y)
\ee
for vector fields $X, \, Y$ in $Q$.
(To avoid confusion with our notation $Q$ for the contact manifold and to highlight that $P$
is not the same tensor field as $Q$ but is  the contact analog thereof,
we use $P$ instead for its notation.)
From \eqref{eq:Q-vanishing}, we have now
$$
 \Pi P(\Pi Z_1,\Pi Z_2) = \frac{1}{2}J((\nabla^{LC}_{\Pi Z_1}J)\Pi Z_2).
$$
According to the remark made in the beginning of the section,
we choose $B_1$ to be
\be\label{eq:B1}
B_1(Z_1,Z_2) = - \Pi P(\Pi Z_1,  \Pi Z_2)=-\frac{1}{2}J((\nabla^{LC}_{\Pi Z_1}J)\Pi Z_2).
\ee

First we consider the induced vector bundle connection
on the Hermitian bundle $\xi \rightarrow Q$, which we denote by $\nabla^{tmp;1,\pi}$: it is defined by
\be\label{eq:nablatmp1pi}
\nabla^{tmp;1,\pi}_X Y: = \pi \nabla^{tmp;1}_X Y
\ee
for a vector field $Y$ tangent to $\xi$, i.e., a section of $\xi$ for arbitrary vector field $X$ on $Q$.
We now prove the $J$ linearity of $\nabla^{tmp;1,\pi}$.
\begin{lemma}\label{lem:J-linear-tmp1} Let $\pi: TQ \to \xi$ be the projection. Then
$
\nabla^{tmp;1,\pi}_{X}(JY)= J\nabla^{tmp;1,\pi}_{X}Y
$
for $Y \in \xi$ and all $X\in TQ$.
\end{lemma}
\begin{proof} For $X\in \xi$,
\bea
\nabla^{tmp;1}_{X}(JY)& = & \nabla^{LC}_{X}(JY) - \Pi P(X, JY)\nonumber\\
&=& (J\nabla^{LC}_{X}Y+(\nabla^{LC}_{X}J)Y)-\frac{1}{2} J((\nabla^{LC}_{X}J)JY)\nonumber\\
&=&J\nabla^{LC}_{X}Y+ (\nabla^{LC}_{X}J)Y
-\frac{1}{2}\Pi((\nabla^{LC}_{X}J)Y)+\frac{1}{2}J((\nabla^{LC}_{X}\Pi)Y)\label{eq:J-linear1}\\
&=& J\nabla^{LC}_{X}Y+ (\nabla^{LC}_{X}J)Y-\frac{1}{2}\Pi((\nabla^{LC}_{X}J)Y)\nonumber
\eea
where we use \eqref{eq:J1} to get the last two terms in the third equality and use \eqref{eq:J2} to see that the last term in
\eqref{eq:J-linear1} vanishes. Hence,
$$
\pi \nabla^{tmp;1}_{X}(JY) = \pi\nabla^{tmp;1}_{X}(JY)
= J\nabla^{LC}_{X}Y+\frac{1}{2}\pi ((\nabla^{LC}_{X}J)Y).
$$
On the other hand, we compute
$$
J\pi \nabla^{tmp;1}_{X}Y = J\left( \nabla^{LC}_{X}Y-\frac{1}{2}J((\nabla^{LC}_{X}J)Y)\right)
= J\nabla^{LC}_{X}Y+\frac{1}{2}\pi ((\nabla^{LC}_{X}J)Y).
$$
Hence we have now $\pi \nabla^{tmp;1}_{X}(JY)=J\pi \nabla^{tmp;1}_{X}Y$
for $X, \ Y\in \xi$.

On the other hand, we notice that $\nabla^{tmp;1}_{X_\lambda}Y=\nabla^{LC}_{X_\lambda}Y$. By using
Proposition \ref{prop:LCXJ}, the equality
$
\pi \nabla^{tmp;1}_{X}(JY)=J\pi \nabla^{tmp;1}_{X}Y
$
also holds for $X=X_\lambda$, and we are done with the proof.
\end{proof}

Next we study the metric property of $\nabla^{tmp;1}$ by computing
$\langle \nabla^{tmp;1}_X Y, Z\rangle + \langle Y, \nabla^{tmp;1}_X Z \rangle$ for arbitrary $X, Y, Z\in TQ$.

Using the metric property of the Levi-Civita connection, we derive
\bea
&{}&\langle \nabla^{tmp;1}_X Y, Z\rangle + \langle Y, \nabla^{tmp;1}_X Z \rangle-X\langle Y, Z\rangle\nonumber\\
&=&\langle \nabla^{LC}_X Y, Z\rangle + \langle Y, \nabla^{LC}_X Z \rangle-X\langle Y, Z\rangle
-\langle \Pi P(\Pi X, \Pi Y), Z\rangle - \langle Y, \Pi P(\Pi X, \Pi Z)\rangle\nonumber\\
&=& -\langle \Pi P(\Pi X, \Pi Y), Z\rangle - \langle Y, \Pi P(\Pi X, \Pi Z)\rangle\label{eq:metric-equality},
\eea
The following lemma shows that when $X, Y, Z\in \xi$ this last line
vanishes. This is the contact analog to Proposition \ref{prop:N=8T} whose proof is
also similar thereto this time based on Lemma \ref{lem:N-property-contact}. Since we
work in the contact case for which we cannot directly quote its proof here, we give
complete proof for readers' convenience.

\begin{lemma}\label{lem:metric-xi}
For $X, Y, Z\in \xi$,
$
\langle P(X, Y), Z\rangle + \langle Y, P(X, Z)\rangle=0.
$
In particular,
\beastar
\langle \nabla^{tmp;1}_X Y, Z\rangle + \langle Y, \nabla^{tmp;1}_X Z \rangle=X\langle Y, Z\rangle.
\eeastar
\end{lemma}
\begin{proof}
We compute for $X, Y, Z\in \xi$,
\bea
&{}&\langle P(X, Y), Z\rangle + \langle Y, P(X, Z)\rangle\nonumber\\
&=&\frac{1}{2}\langle J((\nabla^{LC}_XJ)Y), Z\rangle +\frac{1}{2}\langle Y, J((\nabla^{LC}_XJ)Z)\rangle\nonumber\\
&=&-\frac{1}{2}\langle (\nabla^{LC}_XJ)Y, JZ\rangle -\frac{1}{2}\langle JY, (\nabla^{LC}_XJ)Z\rangle\nonumber\\
&=&-\frac{1}{4}\langle N(Y, JZ), JX\rangle -\frac{1}{4}\langle N(Z, JY), JX\rangle\label{eq:metric-xi1}\\
&=&-\frac{1}{4}\langle \Pi N(Y, JZ)+\Pi N(Z, JY), JX\rangle=0\label{eq:metric-xi2},
\eea
where we use the third equality of Corollary \ref{cor:nablaLCJ} for \eqref{eq:metric-xi1} and use
the second equality of Lemma \ref{lem:N-property-contact} for the vanishing of \eqref{eq:metric-xi2}.
\end{proof}

Now, we are ready to state the following proposition.
\begin{proposition}\label{prop:hermitian-tmp1}
The vector bundle connection $\nabla^{tmp;1,\pi}: = \pi \nabla^{tmp;1}$ is an Hermitian connection of the Hermitian bundle $\xi\rightarrow Q$.
\end{proposition}
\begin{proof}
What is now left to show is that for any $Y, \ Z \in \xi$,
\beastar
\langle \nabla^{tmp;1}_{X_\lambda} Y, Z\rangle + \langle Y, \nabla^{tmp;1}_{X_\lambda} Z \rangle=X_\lambda \langle Y, Z\rangle,
\eeastar
which immediately follows from our construction of $\nabla^{tmp;1}$ since
$$
\nabla^{tmp;1}_{X_\lambda} Y = \nabla^{LC}_{X_\lambda} Y,\quad
\nabla^{tmp;1}_{X_\lambda} Z = \nabla^{LC}_{X_\lambda} Z.
$$
\end{proof}
With direct calculation, one can check 
the metric property when the Reeb direction gets involved.
\begin{lemma}\label{lem:metric-X}
For $Y, \ Z\in \xi$,
$
\langle \nabla^{tmp;1}_Y {X_\lambda}, Z\rangle + \langle X_\lambda, \nabla^{tmp;1}_Y Z \rangle
= 0$.
\end{lemma}

Now we study the torsion property of $\nabla^{tmp;1}$.
Denote the torsion of $\nabla^{tmp;1}$ by $T^{tmp;1}$.
Similar as for the almost Hermitian case, define
$\Theta^\pi=\Pi' T^{tmp;1, \pi}_\C$.
Here we decompose
$$
T^{tmp;1}|_\xi = \pi T^{tmp;1}|_\xi + \lambda(T^{tmp;1, \pi}|_\xi)\, X_\lambda
$$
and denote $T^{tmp;1, \pi}|_\xi:=\pi T^{tmp;1, \pi}|_\xi$,
The proof of the following lemma follows essentially the same strategy as that of
the proof of \cite[Theorem 3.4]{kobayashi-nomizu}.
We give the complete proof for readers' convenience.

\begin{lemma}\label{lem:torsion}
For $Y\in \xi$,
$
T^{tmp;1}(X_\lambda, Y)=0,
$
and
$$
T^{tmp;1, \pi}|_\xi =\frac{1}{4} N^\pi|_\xi, \quad
\lambda(T^{tmp;1}|_\xi)=0.
$$
In particular, $\Theta^\pi|_\xi$ is of $(0,2)$ form.
\end{lemma}

\begin{proof}
Since $\nabla^{tmp;1} = \nabla^{LC} - \Pi P(\Pi, \Pi)$ and $\nabla^{LC}$ is torsion free, we derive for $Y, \ Z \in \xi$,
\beastar
T^{tmp;1}(Y, Z) &= &T^{LC}(Y,Z) -\Pi P(Y,Z)+\Pi P(Z, Y)\\
& = & \frac{1}2J \nabla^{LC}_YJ Z - \frac{1}2J \nabla^{LC}_ZJ Y.
\eeastar
from the general torsion formula.

Next we calculate $-\Pi P(\Pi Y, \Pi Z)+\Pi P(\Pi Z, \Pi Y)$ using the formula
\beastar
\frac{1}2J \nabla^{LC}_YJ Z - \frac{1}2J \nabla^{LC}_ZJ Y &=& \frac{1}{4}\pi([JY,J Z]-\pi [Y,Z]-J[JY, Z]-J[Y, JZ])\\
&=& \frac{1}{4}\pi N(Y,Z).
\eeastar
This follows from the general formula
\be\label{eq:PYZPZY}
-P(Y,Z) + P(Z,Y) = \frac{1}{4}([JY,J Z]-\Pi [Y,Z]-J[JY, Z]-J[Y, JZ]),
\ee
whose derivation we refer \cite[Appendix]{oh-wang}.

On the other hand, since the added terms to $\nabla^{LC}$ only involves $\xi$-directions, the
$X_\lambda$-component of the torsion does not change and so
$$
\lambda(T^{tmp;1}|_\xi) = \lambda(T^{LC}|_\xi) = 0.
$$
This finishes the proof.
\end{proof}

From the definition of $\nabla^{tmp;1}$, we have the following lemma from the properties of the Levi-Civita connection in Proposition \ref{lem:LC-XX}.
\begin{lemma}\label{lem:Xlambda-condition}
$\nabla^{tmp;1}_{X_\lambda}X_\lambda=0$ and
 $\nabla^{tmp;1}_Y X_\lambda\in \xi$ for any $Y\in \xi$.
\end{lemma}
We also get the following property by using Lemma \ref{lem:YXlambda-LC} for Levi-Civita connection.

\begin{lemma}\label{lem:YXlambda-tmp1}
For any $Y\in \xi$, we have
$\nabla^{tmp;1}_Y X_\lambda=\frac{1}{2}JY+\frac{1}{2}(\CL_{X_\lambda}J)JY$.
\end{lemma}

We end the construction of $\nabla^{tmp;1}$ by summarizing that $\nabla^{tmp;1}$ satisfies Axioms (1)-(4),(6) and (5;$-1$), i.e., $\nabla^{tmp;1} = \nabla^{\lambda;-1}$.

\subsection{Modification 2; $\nabla^{tmp;2}$}
\label{subsec:modify2}
\par

Now we introduce another modification $\nabla^{tmp;2}$ starting from $\nabla^{tmp;1}$ to make it satisfy Axiom
(5;$c$) and preserve other axioms for any given constant $c \in \R$.
We define $\nabla^{tmp;2} = \nabla^{tmp;1} + B_2$ for the tensor $B_2$ given as
\be\label{eq:B2}
B_2(Z_1,Z_2) = \frac{1}{2}(1+c)\left(
-\left\langle Z_2, X_\lambda\right\rangle JZ_1
-\left\langle Z_1, X_\lambda\right\rangle JZ_2
+ \left\langle JZ_1, Z_2\right\rangle X_\lambda\right).
\ee


\begin{proposition} The connection $\nabla^{tmp;2}$ satisfies all the properties of
the canonical connection with constant $c$. In particular $\nabla: = \nabla^{tmp;2}$ with $c = 0$
is the contact triad connection.
\end{proposition}
\begin{proof}





The checking of all Axioms are straightforward, and we only do it for Axiom (5;$c$) here.
\beastar
&{}&\nabla^{tmp;2}_{JY} X_\lambda+J\nabla^{tmp;2}_Y X_\lambda\\
&=&\nabla^{tmp;1}_{JY} X_\lambda-\frac{1}{2}\left(1+c\right)JJY
+J\nabla^{tmp;1}_Y X_\lambda-J\frac{1}{2}\left(1+c\right)JY\\
&=&-Y+(1+c)Y=cY.
\eeastar

\end{proof}




Before ending this section, we restate the following properties which will be useful for calculations
involving contact Cauchy-Riemann maps performed in \cite{oh-wang1}, \cite{oh-wang2}.

\begin{proposition}\label{prop:YX-triad-connection} Let $\nabla$ be the connection
satisfying Axiom (1)-(4),(6) and (5; $c$), then
$
\nabla_Y X_\lambda=-\frac{1}{2}cJY+\frac{1}{2}(\CL_{X_\lambda}J)JY
$.
In particular, for the contact triad connection,
$
\nabla_Y X_\lambda=\frac{1}{2}(\CL_{X_\lambda}J)JY
$.
\end{proposition}
\begin{proof}
We already gave its proof in the last part of Section \ref{sec:can-connection}.
\end{proof}

\begin{proposition}\label{eq:torsion3} Decompose the torsion of
$\nabla$ into $T = \pi T + \lambda(T)\, X_\lambda$.
The triad connection $\nabla$ has its torsion given by $T(X_\lambda, Z) = 0$ for all $Z \in TQ$, and
for all $Y, \, Z \in \xi$,
\beastar
\pi T(Y,Z) & = & \frac{1}{4}\pi N(Y,Z) = \frac{1}{4}\left((\CL_{JY} J)Z + (\CL_Y J)JZ\right) \\
\lambda(T(Y,Z)) & = & d\lambda(Y, Z).
\eeastar
\end{proposition}
\begin{proof} We have seen
$
\pi T^{tmp;2}|_\xi =  \pi T^{tmp;1}|_\xi = \frac{1}{4} N^\pi|_\xi
$.
On the other hand, a simple computation shows
$
N^\pi(Y,Z) = (\CL_{JY}J)Z - J(\CL_Y J)Z = (\CL_{JY} J)Z + (\CL_Y J)JZ
$, which proves the first equality.

For the second, a straightforward computation shows
\beastar
\lambda(T^{tmp;2}(Y,Z))  =  \lambda(T^{tmp;1}(Y,Z)) + (1+c)\, \langle JY, Z \rangle
 =  (1+ c)\, d\lambda(Y,Z)
\eeastar
for general $c$. Substituting $c = 0$, we obtain the second equality.
This finishes the proof.
\end{proof}

\begin{acknowledgement}
We thank Luigi Vezzoni and Liviu Nicolaescu for alerting their
works \cite{vezzoni} and \cite{nicola} respectively on special connections on contact manifolds after the original version of the present paper was posted in the arXiv e-print.
We also thank them for helpful discussions. 
The present work was supported by the IBS project IBS-R003-D1.
\end{acknowledgement}
%
%
%


\end{document}